\documentclass[11pt]{article}

\usepackage{amscd,amsmath, amssymb, fancyhdr, url}

\usepackage{fourier, libertine}

\newcommand{\version}{version 3.0,\ \ March 03, 2019}
\numberwithin{equation}{section}

\def\eqref#1{(\ref{#1})}

\newcommand{\Z}{{\mathbb Z}}
\newcommand{\C}{{\mathbb C}}
\newcommand{\R}{{\mathbb R}}

\def\1{\sqrt{-1}\:}

\newcommand{\cntrct}                % contraction with a vector field
{\hspace{2pt}\raisebox{1pt}{\text{$\lrcorner$}}\hspace{2pt}}

\renewcommand{\bar}{\overline}
\renewcommand{\phi}{\varphi}
\renewcommand{\epsilon}{\varepsilon}
\renewcommand{\geq}{\geqslant}

\newcommand{\const}{\operatorname{\text{\sf const}}}

\newcommand{\Aut}{\operatorname{Aut}}

\newcommand{\codim}{\operatorname{codim}}

\newcounter{Mycounter}[section]
\newcounter{lemma}[section]
\renewcommand{\thelemma}{{Lemma \thesection.\arabic{lemma}}}
\newcommand{\lemma}{%
     \setcounter{lemma}{\value{Mycounter}}
     \refstepcounter{lemma}
     \stepcounter{Mycounter}
     {\noindent \bf \thelemma:\ }}

\newcounter{claim}[section]
\renewcommand{\theclaim}{{Claim \thesection.\arabic{claim}}}
\newcommand{\claim}{%
     \setcounter{claim}{\value{Mycounter}}
     \refstepcounter{claim}
     \stepcounter{Mycounter}
     {\noindent \bf \theclaim:\ }}

\newcounter{sublemma}[section]
\newcounter{corollary}[section]
\renewcommand{\thecorollary}{{Corollary \thesection.\arabic{corollary}}}
\newcommand{\corollary}{%
     \setcounter{corollary}{\value{Mycounter}}
     \refstepcounter{corollary}
     \stepcounter{Mycounter}
     {\noindent \bf \thecorollary:\ }}

\newcounter{theorem}[section]
\renewcommand{\thetheorem}{{Theorem \thesection.\arabic{theorem}}}
\newcommand{\theorem}{%
     \setcounter{theorem}{\value{Mycounter}}
     \refstepcounter{theorem}
     \stepcounter{Mycounter}
     {\noindent \bf \thetheorem:\ }}

\newcounter{conjecture}[section]
\newcounter{proposition}[section]
\renewcommand{\theproposition} {{Proposition \thesection.\arabic{proposition}}}
\newcommand{\proposition}{%
     \setcounter{proposition}{\value{Mycounter}}
     \refstepcounter{proposition}
     \stepcounter{Mycounter}
     {\noindent \bf \theproposition:\ }}

\newcounter{definition}[section]
\renewcommand{\thedefinition} {{Definition~\thesection.\arabic{definition}}}
\newcommand{\definition}{%
     \setcounter{definition}{\value{Mycounter}}
     \refstepcounter{definition}
     \stepcounter{Mycounter}
     {\noindent \bf \thedefinition:\ }}

\newcounter{example}[section]
\newcounter{remark}[section]
\renewcommand{\theremark}{{Remark \thesection.\arabic{remark}}}
\newcommand{\remark}{%
     \setcounter{remark}{\value{Mycounter}}
     \refstepcounter{remark}
     \stepcounter{Mycounter}
     {\noindent \bf \theremark:\ }}

\newcounter{problem}[section]
\newcounter{question}[section]
\makeatletter

\@addtoreset{equation}{section}
\@addtoreset{footnote}{section}
\makeatother

\def\blacksquare{\hbox{\vrule width 5pt height 5pt depth 0pt}}
\def\endproof{\blacksquare}

\begin{document}

\begin{center}
{\LARGE\bf
Hopf surfaces in locally conformally K\"ahler manifolds with potential}\\[3mm]
%%%%%%%%%%%%%%%%%%%%%%%%%%%%%%%%%%%%%%%%%%%%%%%%%%%%%%%%%%%%
{\large
Liviu Ornea\footnote{Partially supported by CNCS UEFISCDI, project number PN-II-ID-PCE-2011-3-0118.}, and 
Misha
Verbitsky\footnote{Partially supported 
by  the  Russian Academic Excellence Project '5-100',
FAPERJ E-26/202.912/2018 and CNPq - Process 313608/2017-2.\\[1mm]

}\\[4mm]
}

\end{center}

\hfill

\hfill

%%%%%%%%%%%%%%%%%%%%%%%%%%%%%%%%%%%%%%%%%%%%%%%%
{\small
\hspace{0.15\linewidth}
\begin{minipage}[t]{0.7\linewidth}
{\bf Abstract} \\ 
An LCK manifold with potential is a quotient $M$ of a
K\"ahler manifold $X$ equipped with a positive plurisubharmonic function 
$f$, such that the monodromy group acts on $X$ by holomorphic homotheties
and maps $f$ to a function proportional to $f$.
It is known that a compact $M$ admits an LCK potential
if and only if it can be holomorphically embedded to a Hopf manifold.
We prove that any non-Vaisman, compact LCK manifold with potential
contains a complex surface (possibly singular)
with normalization biholomorphic to a Hopf surface 
$H$. Moreover, $H$ can be chosen non-diagonal, 
hence, also not admitting a Vaisman 
structure. 
\end{minipage}
}
%%%%%%%%%%%%%%%%%%%%%%%%%%%%%%%%%%%%%%%%%%%%%%%%

\tableofcontents

%%%%%%%%%%%%%%%%%%%%%%%%%%%%%%%%%%%%%%%%%%%%%%%%%%%%%%%%%%%%%%%%%%%%%%%%

\section{Introduction: LCK manifolds}

%%%%%%%%%%%%%%%%%%%%%%%%%%%%%%%%%%%%%%%%%%%%%%%%%%%%%%%%%%%%%%%%%%%%%%%%

Let $(M,I)$ be a complex manifold, $\dim_\C M\geq 2$. It
is called {\bf locally conformally K\"ahler} (LCK) if it
admits a  Hermitian metric $g$ whose 
fundamental 2-form $\omega(\cdot,\cdot):=g(\cdot, I\cdot)$
satisfies
\begin{equation}\label{deflck}
d\omega=\theta\wedge\omega,\quad d\theta=0,
\end{equation}
for a certain closed 1-form $\theta$ called {\bf the Lee
form}. 

Definition \eqref{deflck} is equivalent to the existence
of a  covering $\tilde M$ endowed with a K\"ahler metric $\Omega$ which is
acted on by the deck group $\Aut_M(\tilde M)$ by homotheties. Let 
\begin{equation}\label{chi}
\chi:\Aut_M(\tilde M)\longrightarrow \R^{>0},\quad \chi(\tau)=\frac{\tau^*\Omega}{\Omega},
\end{equation}
be the group homomorphism which associates to a homothety its scale factor.

For definitions and examples, see \cite{DO} and our more recent papers.

An LCK manifold $(M,\omega, \theta)$ is called {\bf
  Vaisman} if $\nabla\theta=0$, where $\nabla$ is the
Levi-Civita connection of $g$. The main example of Vaisman manifold is the
diagonal Hopf manifold (\cite{_OV:Shells_}). The Vaisman compact complex
surfaces are classified in \cite{be}. 

Note that there exist compact LCK manifolds which do not admit Vaisman metrics. Such are the LCK Inoue surfaces, \cite{be}, the Oeljeklaus-Toma manifolds, \cite{kas}, and the non-diagonal Hopf manifolds, \cite{_OV:Shells_}.

  It is known that on any Vaisman manifold with Lee form normalized to have length 1, the following formula holds, \cite{vai}, \cite{DO}:
\begin{equation}\label{dc}
d\theta^c=\theta\wedge\theta^c-\omega, \quad \text{where}\quad \theta^c(X)=-\theta(IX).
\end{equation}
Moreover, one can see, \cite{_Verbitsky:LCHK_}, that the
(1,1)-form $\omega_0:=-d^c\theta$ is semi-positive
definite, having all eigenvalues\footnote{The eigenvalues
  of a Hermitian form $\eta$  are the eigenvalues of the
  symmetric operator $L_\eta$ defined by the equation
  $\eta(x,Iy)=g(Lx,y)$.} 
positive, but one which is 0.

An LCK manifold is called {\bf with potential} if it
admits a K\"ahler covering on which the K\"ahler metric
has a global and positive  potential function which
is acted on by holomorphic homotheties by the deck group. Among the
examples:  all Vaisman manifolds, but also non-Vaisman
ones, such as the non-diagonal Hopf manifolds, \cite{ov01}, \cite{ov_rank}. On the other hand, there exist  compact LCK manifolds which cannot admit LCK potential, {\em e.g.} Inoue surfaces (see \cite{ov01}) and their higher dimensional analogues, the Oeljeklaus-Toma manifolds, see \ref{nopot}.

One can prove, \cite{ov_rank}, that on a compact manifold,
a positive, automorphic potential can always be deformed
to a {\bf proper} positive, automorphic potential. The
existence of such a potential is equivalent with the image
of the character $\chi$ being isomorphic with $\Z$. In
this case, the LCK manifold with potential is called {\bf
  of LCK rank 1}.

On the K\"ahler covering of an LCK manifold with
potential, one has $\pi^*\omega=\psi^{-1}dd^c\psi$, where
the potential is $\psi=e^{-\nu}$ and the Lee form is
$\pi^*\theta=d\nu$. Hence we have  (\cite{ov_imrn_10},
also \cite{ad}):

\hfill

\claim \label{dc_on_pot}
Equation \eqref{dc} is satisfied on LCK manifolds with potential.

\hfill

The aim of this paper is to prove that any compact non-Vaisman LCK
manifold with potential contains  a complex surface (possibly singular)
with normalization biholomorphic to a Hopf surface 
which is non-diagonal (this is the same as being non-Vaisman). 
As an application, we show that a compact LCK manifold with
potential $(M,\omega, \theta)$ is Vaisman if and only if
the form $d^c\theta$ is sign semi-definite.

%%%%%%%%%%%%%%%%%%%%%%%%%%%%%%%%%%%%%%%%%%%%%%%%%%%%%%%%%%%%%%
\section{The form $\omega_0$ on a compact LCK manifold with potential}
%%%%%%%%%%%%%%%%%%%%%%%%%%%%%%%%%%%%%%%%%%%%%%%%%%%%%%%%%%%%%

In general, on an LCK manifold $M$ with potential $\psi$ on $\tilde M$, the
norm (w.r.t. the LCK metric) of the Lee form $d\psi$ is
not constant. The constancy of the norm of the Lee form 
is equivalent to the LCK metric being Gauduchon (see 
\ref{_Gauduchon_Lee_length_Proposition_} below) and
Vaisman, as shown in \cite{mm}.

\hfill

%%%%%%%%%%%%%%%%%%%%%%%%%%%%%%%%%%%%%%%%%%%%%%%%%%%%%%%%%%%%
\definition
On a complex manifold of complex dimension $n$, a
Hermitian metric whose Hermitian 2-form $\omega$ satisfies the equation
$\partial\bar\partial\omega^{n-1}=0$ is called {\bf Gauduchon}.

\hfill

\remark 
On a compact Hermitian manifold, a Gauduchon metric exists
in each conformal class and it is unique up to
homothety. Moreover, it is characterized by the
co-closedness of its Lee form. A Vaisman metric {\em is}
a Gauduchon metric in its conformal class, \cite{gau}.

\hfill

%%%%%%%%%%%%%%%%%%%%%%%%%%%%%%%%%%%%%%%%%%%%%%%%%%%%%%%%%%%%
\proposition\label{_Gauduchon_Lee_length_Proposition_}
Let $(M, \omega, \theta)$ be a compact LCK manifold with   
potential. 
Then the LCK form $\omega$ 
is Gauduchon if and only if $|\theta|=\const$.

\hfill

 {\bf Proof:} The Hermitian form $\omega$ is Gauduchon if and only if 
$dd^c\omega^{n-1}=0$.

We compute $dd^c\omega^{n-1}$ using
equation \eqref{dc} which is satisfied on an LCK manifold
with potential (\ref{dc_on_pot}). This gives
\[
dd^c\omega^{n-1}=(n-1)^2\omega^{n-1}\wedge\theta\wedge\theta^c+(n-1)\omega^{n-1}\wedge
d\theta^c.
\]
On the other hand,
$$\omega^{n-1}\wedge\theta\wedge\theta^c=\frac 1n |\theta|^2\omega^n$$
and
$$d\theta^c\wedge\omega^{n-1}=-\omega\wedge\omega^{n-1}+
\theta\wedge\theta^c\wedge\omega^{n-1}=\left(\frac 1n |\theta|^2-1\right)\omega^n.$$
All in all we get:
\[ dd^c\omega^{n-1}=\frac{(n-1)^2}{n}|\theta|^2\omega^n+(n-1)\left(\frac
1n|\theta|^2-1\right)\omega^n=(n-1)\big(|\theta|^2-1\big)\omega^n.
\]
Then $dd^c\omega^{n-1}=0$ if and only if $|\theta|=1$. This finishes
the proof.   \endproof

\hfill

Observe now that the eigenvalues of $\omega_0=-d\theta^c$
are $1$ (with multiplicity $n-1$) and $1-|\theta|^2$. As
$\omega_0$ is exact on a compact manifold, its top power
cannot be sign-definite (Stokes theorem). Two
possibilities occur:
\begin{enumerate}
\item $|\theta|$ is non-constant, and then  $1-|\theta|^2$ has to change sign on $M$;
\item $|\theta|=\const.$ and then $|\theta|=1$ and $\omega_0$ is semi-positive definite.
\end{enumerate}

  We obtained the following corollary.

\hfill

\corollary\label{pot_gau}
Let $(M,\omega,\theta)$ be a compact LCK manifold with
potential. Then the LCK metric is Gauduchon if and only if
$\omega_0=-d\theta^c$ is semi-positive definite, and is then Vaisman.

\hfill

\remark 
Our interest in studying the form $\omega_0$ on
compact LCK manifolds arose from the attempt to clarify
the relation between the pluricanonical condition
($(\nabla\theta)^{1,1}=0$, equivalent with 
$d^c\theta=\theta\wedge\theta^c-|\theta|^2\omega$),
introduced in \cite{ko}, and the existence of a positive,
automorphic potential. In fact, from the above it can be
seen that a compact LCK manifold with potential, and with
constant norm of $\theta$ is pluricanonical, and hence,
by \cite{mm}, it is Vaisman.

%%%%%%%%%%%%%%%%%%%%%%%%%%%%%%%%%%%%%%%%%%%%%%%%%%%%%%%%%%%%%%%%%%%%%%%%

\section{Hopf surfaces in LCK manifolds with potential}

%%%%%%%%%%%%%%%%%%%%%%%%%%%%%%%%%%%%%%%%%%%%%%%%%%%%%%%%%%%%%%%%%%%%%%%%

%%%%%%%%%%%%%%%%%%%%%%%%%%%%%%%%%%%%%%%%%%%%%%%%%%%%%%%%%%%%%%%%%%%%%%%%
\subsection{Complex surfaces of K\"ahler rank 1}\label{_Surfaces_Section_}
%%%%%%%%%%%%%%%%%%%%%%%%%%%%%%%%%%%%%%%%%%%%%%%%%%%%%%%%%%%%%%%%%%%%%%%%

\definition {\bf (\cite{hl})} A compact complex surface is {\bf of
K\"ahler rank 1} if and only if it is not K\"ahler but it
admits a closed semipositive (1,1)-form whose zero
locus is contained in a curve.

\hfill

\lemma A  compact  LCK surface $M$ with potential and
semi-positive form $\omega_0$  has K\"ahler rank 1.

\hfill

{\bf Proof:} We have to show that $M$ cannot admit a
K\"ahler metric. By absurd, if  $M$ admitted a K\"ahler
form $\Omega$, then, as $\omega_0$ is exact,
$\omega_0\wedge\Omega=-d(\theta^c\wedge\Omega)$ was an
exact volume form, which is impossible by Stokes'
theorem. Hence $M$ is non-K\"ahler.
\endproof

\hfill

Recall that a {\bf Hopf surface}
is a finite quotient of $H$, where $H$ is a quotient of $\C^2 \backslash 0$
by a polynomial contraction. 
A Hopf surface is {\bf diagonal} if this 
polynomial contraction is expressed by a 
diagonal matrix. 

\hfill

Compact surfaces of K\"ahler rank 1 have been 
classified in \cite{ct} and \cite{b}. They can be:
\begin{enumerate}
\item Non-K\"ahler elliptic fibrations, 
\item Diagonal Hopf surfaces and their blow-ups, 
\item Inoue surfaces and their blow-ups.
\end{enumerate}

The LCK Inoue surfaces cannot have  LCK metrics with potential, as shown in \cite[Corollary 3.13]{aot}.

A cover of a blow-up of any complex manifold cannot admit plurisubharmonic functions because, by the lifting criterion, the projective spaces  contained in the blow-up lift to the cover. Thus blow-ups cannot have global potential.

We are left with non-K\"ahler elliptic fibrations and diagonal Hopf surfaces which are known to admit Vaisman metrics, see {\em e.g.} \cite{be}. And hence:

\hfill

\proposition\label{surfvai} 
All compact  LCK surfaces with potential and with
semi-positive  form $\omega_0$ are  Vaisman.
\endproof

\hfill

For further use it is convenient to list 
all criteria used to distinguish Vaisman
Hopf surfaces from non-Vaisman ones.

\hfill

%%%%%%%%%%%%%%%%%%%%%%%%%%%%%%%%%%%%%%%%%%%%%%%%%%%%%%%%%%%%
\theorem\label{_Hopf_surface_Vaisman_Theorem_}
Let $M$ be a Hopf surface. Then the following
are equivalent.
\begin{description}
\item[(i)] $M$ is Vaisman.
\item[(ii)] $M$ is diagonalizable.
\item[(iii)] $M$ has K\"ahler rank 1.
\item[(iv)] $M$ contains at least two distinct elliptic curves.
\end{description}

{\bf Proof:} The equivalence of the first three conditions is proven above.
The equivalence of (iv) and (ii) is shown by Iku Nakamura and Masahide Kato 
(\cite[Theorem 5.2]{_Nakamura:curves_}). Note that the cited result refers to primary Hopf surfaces, but we can always pass to a finite covering
and  the number of elliptic curves will not change
 because the eigenvectors for rationally independent eigenvalues  cannot be mutually exchanged,
 and if they were dependent, they would produce infinitely many elliptic curves. \endproof

%%%%%%%%%%%%%%%%%%%%%%%%%%%%%%%%%%%%%%%%%%%%%%%%%%%%%%%%%%%%%%%%%%%%%%%%

%%%%%%%%%%%%%%%%%%%%%%%%%%%%%%%%%%%%%%%%%%%%%%%%%%%%%%%%%%%%%%%%%%%%%%%%
\subsection{Algebraic groups and the Jordan-Chevalley decomposition}
%%%%%%%%%%%%%%%%%%%%%%%%%%%%%%%%%%%%%%%%%%%%%%%%%%%%%%%%%%%%%%%%%%%%%%%%

In this section we fix an $n$-dimensional complex vector
space $V$.

\hfill

%%%%%%%%%%%%%%%%%%%%%%%%%%%%%%%%%%%%%%%%%%%%%%%%%%%%%%%%%%%%
\lemma 
Let $A\in \mathrm{GL}(V)$ be a linear operator,  and
$\langle A \rangle$ the group generated by $A$. Denote by
$G$ the Zariski closure of $\langle A \rangle$ in $\mathrm{GL}(V)$.
Then, for any $v\in V$, the Zariski closure $Z_v$
of the orbit $\langle A\rangle\cdot v$ is equal to the usual
closure of $G\cdot v$.

\hfill

{\bf Proof:} Clearly, $Z_v$ is $G$-invariant. Indeed,
its normalizer $N(Z_v)$ in $GL(V)$ is an algebraic group containing $\langle A
\rangle$, hence $N(Z_v)$ contains $G$. The converse is also
true: since  $\langle A\rangle$ normalizes  $\langle A\rangle\cdot v$,
its Zariski closure $G$ normalizes the Zariski closure
$Z_v$ of the orbit. Therefore, the orbit  $G\cdot v$
is contained in $Z_v$. Since $G\cdot v$ is a constructible
set, its Zariski closure coincides with its usual closure, \cite{hum}, \cite{kol}.
This gives $\overline{G\cdot v}\subset Z_v$. As $\overline{G\cdot v}$ is
algebraic and contains $\langle A\rangle\cdot v$, the inclusion $Z_v\subset \overline{G\cdot v}$ is also true.
\endproof

\hfill

The reason we take the Zariski closure is explained in the following (see also \cite[Theorem 2.1]{ov_imrn_10}):

\hfill

\claim Let $I\subset \C[z_1, ..., z_n]$
be an ideal which is invariant with respect
to  an isomorphism $A$ of the space
$\langle z_1, \ldots, z_n\rangle$ acting
on the polynomial ring. Then $I$
is invariant with respect to the
Zariski closure $G$ of $\langle A\rangle$.

\hfill

{\bf Proof:} First, we show that the 0-adic completion $\hat I$ of
$I$ is $G$-invariant in the 0-adic completion
of the polynomial ring $\C[z_1, ..., z_n]$, which is the
ring of formal power series $\C[[z_1, ..., z_n]]$. However, any
$A$-invariant subspace in a finite-dimensional
space is $G$-invariant by definition of $G$, and
the ideal $\hat I$ is obtained as an inverse
limit of finite-dimensional subspaces of
finite quotients of the polynomial ring.
Therefore, $\hat I$ is $G$-invariant.
The ideal $I$ is $G_A$-invariant, 
because
$I=\hat I \cap \C[z_1, ..., z_n]$. \endproof

\hfill 

Let now $G\subset \mathrm{GL}(V)$ be an algebraic group over $\C$.
Recall that an element $g\in G$ is called {\bf semisimple}
if it is diagonalizable, and {\bf unipotent} if $g = e^n$,
where $n$ is a nilpotent element of its Lie algebra.

\hfill

%%%%%%%%%%%%%%%%%%%%%%%%%%%%%%%%%%%%%%%%%%%%%%%%%%%%%%%%%%%%
\theorem
(Jordan-Chevalley decomposition, \cite[Section 15]{hum})  
For any algebraic group $G\subset \mathrm{GL}(V)$,  any $g\in G$ can be
represented as a product of two commuting elements
$g=g_sg_u$, where $g_s$ is semisimple, and $g_u$
unipotent. Moreover, this decomposition is unique and
functorial under homomorphisms of algebraic groups.

\hfill

%%%%%%%%%%%%%%%%%%%%%%%%%%%%%%%%%%%%%%%%%%%%%%%%%%%%%%%%%%%%
\corollary\label{_algebraic_action_on_LCK_with_pot_Corollary_}
Let $M$ be a submanifold of a linear Hopf manifold
$H=(V\backslash 0)/A$, $\tilde M \subset V\backslash 0$
its $\Z$-covering, and $G$ the Zariski closure of 
$\langle A \rangle$ in $\mathrm{GL}(V)$. Then $\tilde M$ contains
the $G$-orbit of each point $v\in \tilde M$. Moreover,
$G$ is a product of $G_s:=(\C^*)^k$ and a unipotent group
$G_u$ commuting with $G_s$, and bothe of these groups
preserve $\tilde M$.

\hfill

{\bf Proof:}
Let $X$ be the closure of $\tilde M$ in $\C^N$.
The ideal $I_X$ of $X$ is generated 
by polynomials, as shown in  \cite[Proof of Theorem
  3.3]{ov01}. As the polynomial ring is 
Noetherian,  $I_X$ is finitely generated, \cite{am}.  Therefore, 
$X$ is a cone of a projective variety.

Consider  the smallest algebraic group $G$
containing $A$. Then $G$ acts
naturally on $X$ and preserves it. 
The last assertion of 
\ref{_algebraic_action_on_LCK_with_pot_Corollary_}
is implied by the Jordan-Chevalley decomposition.
\endproof

%%%%%%%%%%%%%%%%%%%%%%%%%%%%%%%%%%%%%%%%%%%%%%%%%%%%%%%%%%%%
\subsection{Finding surfaces in LCK manifolds with potential}
%%%%%%%%%%%%%%%%%%%%%%%%%%%%%%%%%%%%%%%%%%%%%%%%%%%%%%%%%%%%

%%%%%%%%%%%%%%%%%%%%%%%%%%%%%%%%%%%%%%%%%%%%%%%%
\lemma\label{_surface_exists_Lemma_}
Let $M$ be a non-Vaisman submanifold of a linear Hopf manifold
$H=(V\backslash 0)/A$, $\dim_\C M \geq 3$, and $G=G_sG_u$ the Zariski
closure of $\langle A \rangle$ with its Jordan-Chevalley 
decomposition. Then $M$ contains a surface $M_0$, possibly singular, 
with $G_u$ acting non-trivially on its $\Z$-covering 
$\tilde M_0\subset V$.

\hfill

{\bf Proof:} Another form of this statement 
is proven by Masahide Kato (\cite{_Kato:subvarieties_}).

We shall use induction on dimension of $M$.
To prove \ref{_surface_exists_Lemma_} it would suffice
to find a subvariety $M_1 \subset M$ of codimension 1
such that $G_u$ acts non-trivially on its $\Z$-covering 
$\tilde M_1\subset \C^n \backslash 0$ (note that $G_u$ is non-trivial because $M$ is non-Vaisman). Replacing 
$V$ by the smallest $A$-invariant subspace containing
$\tilde M$, we may assume that the intersection
$\tilde M \cap V_1\neq V_1$ for each proper subspace
$V_1\subset V$. Now take a codimension 1 
subspace $V_1\subsetneq V$ which is $A$-invariant and 
 such that $G_u$ acts on $V_1$ non-trivially
(equivalently, such that $A$ acts on $V_1$
non-diagonally). Using the Jordan decomposition
of $A$, such $V_1$ is easy to construct.
Then $\tilde M_1':= V_1 \cap \tilde M$ 
gives a subvariety of $M$ of codimension 1
and with non-trivial action of $G_u$.
\endproof

\hfill

The same argument gives the following corollary,
also parallel to a theorem by Ma. Kato.

\hfill

%%%%%%%%%%%%%%%%%%%%%%%%%%%%%%%%%%%%%%%%%%%%%%%%%%%
\corollary\label{flag}
Let $M$ be a compact LCK manifold with potential.
Then $M$ has a flag of embedded subvarieties
$M\supset M_1\supset M_2\supset \cdots\supset M_{\dim M-1}$
with $\codim M_i =i$.
\endproof

\hfill

%%%%%%%%%%%%%%%%%%%%%%%%%%%%%%%%%%%%%%%%%%%%%%%%%

Recall that Oeljeklaus-Toma manifolds (see \cite{_Oeljeklaus_Toma_}) do not admit complex curves
(\cite{_Verbitskaya:OT-curves_} where the argument doesn't need smoothness). Then \ref{flag} implies (see also \cite{_istrati_otiman_} for a more recent different proof):

\hfill

\corollary\label{nopot} The Oeljeklaus-Toma manifolds
cannot admit LCK structures with potential.
\endproof

\hfill

%%%%%%%%%%%%%%%%%%%%%%%%%%%%%%%%%%%%%%%%%%%%%%%%
\lemma\label{_surface_non-diagonal_Lemma_}
Let $M_1\subset H=(V\backslash 0)/\langle A\rangle$
be a surface in a Hopf manifold, possibly singular,
and $G=G_sG_u$ the Zariski closure of $\langle A \rangle$ 
with its Jordan-Chevalley 
decomposition. Assume that $G_u$ acts on the $\Z$-covering
$\tilde M$ non-trivially. Then the normalization of
$M$ is a non-diagonal Hopf surface.

\hfill

{\bf Proof:}  Replacing $G$ by its quotient
by the subgroup acting trivially on $\tilde M$ if necessary, 
 we may assume that
$G$ acts properly on a general orbit in $\tilde M$. 
Then $G$ is at most 2-dimensional. However, it cannot be
1-dimensional because $G_s$ contains contractions
(hence cannot be 0-dimensional) and $G_u$ acts
non-trivially. Therefore, $G_s\simeq\C^*$ and $G_u\simeq\C$.

Since $G_s$ acts by contractions, the quotient
$S:= \tilde M/G_s$ is a compact curve, equipped with
$G_u$-action which has a dense orbit. The group
$G_u$ can act non-trivially only on a genus 0 curve,
and there is a unique open orbit $O$ of $G_u$, with 
$S \backslash O$ being one point. 

Let now $M$ be a normalization of $M_1$. Since the singular
set of $M_1$ is $G_s$-invariant, it has dimension at least 1, 
and since $M$ is normal, it is non-singular in codimension 1,
hence smooth.

All complex subvarieties
of $M$ are by construction $G$-invariant, and the
complement of an open orbit is an elliptic curve,
hence $M$ has only one elliptic curve. As $M$ is a surface of a Hopf manifold, it is LCK with potential and hence it is a deformation of a Vaisman surface (\cite{ov_imrn_10}) which can be Hopf or elliptic (\cite{be}).
By the classification of the 
non-K\"ahler compact surfaces, a smooth deformation of a non-Hopf 
elliptic surface is again an elliptic surface, and hence it has many elliptic curves. As $M$ has only one elliptic curve, it must be a deformation of a Hopf surface and it is
non-diagonalizable by 
\ref{_Hopf_surface_Vaisman_Theorem_}.
\endproof

\hfill

\theorem \label{surface}
Let $M$ be a non-Vaisman compact LCK manifold with potential, $\dim_\C M\geq 3$. Then $M$ 
contains a surface with normalization biholomorphic to a non-diagonal
Hopf surface.

\hfill

{\bf Proof:}
Let $M$ be a compact LCK manifold with potential, $\dim_\C
M\geq 3$. Then  $M$ is holomorphically embedded
into a Hopf manifold $\C^N\backslash 0 /\langle A\rangle$,
where $A\in \mathrm{GL}(N,\C)$  is a linear operator, see  \cite[Theorem 3.4]{ov01}.
 Applying \ref{_surface_exists_Lemma_} 
and \ref{_surface_non-diagonal_Lemma_},
we find a non-diagonal Hopf surface in $M$.
\endproof

\hfill

As an application, we now prove the following characterization of Vaisman manifolds:

\hfill

\corollary\label{main} 
Let $(M,  I)$ be a compact LCK manifold with
potential. Assume the Hermitian form $\omega_0$ is
semi-positive definite. Then the LCK metric of $(M,I)$ is  Vaisman.

\hfill

{\bf  Proof:} If $\dim_\C M=2$, this is just \ref{surfvai}.
 
If $\dim_\C M\geq 3$, 
$M$ contains a surface whose normalization is a non-diagonal Hopf surface
$H$. Then $\omega_0$ restricts to a semi-positive
definite (1,1) form on $H$. By \ref{surfvai}, $H$ is
Vaisman, and hence diagonal, contradiction. \endproof

\hfill

\noindent{\bf Acknowledgments.} L.O. thanks the Laboratory
for Algebraic Geometry at the Higher School of Economics
in Moscow for hospitality and excellent research
environment during February and April 2014, and April
2015. 

Both authors are indebted to Paul 
Gauduchon, Andrei Moroianu, and Victor Vuletescu 
for extremely useful disussions.

{\small

\noindent {\sc Liviu Ornea\\
{\sc University of Bucharest, Faculty of Mathematics, \\14
Academiei str., 70109 Bucharest, Romania, \emph{and}\\
Institute of Mathematics ``Simion Stoilow'' of the Romanian Academy,\\
21, Calea Grivitei Street
010702-Bucharest, Romania }\\
\tt lornea@fmi.unibuc.ro, \ \ liviu.ornea@imar.ro
}

\hfill

\noindent {\sc Misha Verbitsky\\
{\sc Instituto Nacional de Matem\'atica Pura e
              Aplicada (IMPA) \\ Estrada Dona Castorina, 110,
Jardim Bot\^anico, CEP 22460-320,\\ Rio de Janeiro, RJ - Brasil }\\
also:\\
{\sc Laboratory of Algebraic Geometry, HSE University,\\
Department of Mathematics, 9 Usacheva str., Moscow, Russia,}\\
\tt  verbit@mccme.ru
}}

\end{document}